\documentclass[11pt,bezier]{article}
\usepackage{amsmath,amssymb,amsfonts,euscript,graphicx}

\title{\LARGE\bf   Commutative Local Rings whose  Ideals   are Direct Sums of
Cyclic Modules\thanks {The research
 of the first author was in part supported by
a grant from IPM (No. 91130413). This research is partially
carried out in the IPM-Isfahan Branch. }
\thanks {2010{ \it Mathematics Subject Classification}. Primary 13C05, 13E05, 13F10,
  Secondary  13E10, 13H99. }}

\author{{\bf M. Behboodi$^{{\rm a,b}}$\thanks{Corresponding author.}~ and~} {\bf S. H. Shojaee$^{{\rm a}}$}\\
 {\small{ $^{{\rm a}}$Department of Mathematical Sciences, Isfahan University of Technology}}\vspace{-1mm}\\
  {\small{ P.O.Box: 84156-83111, Isfahan, Iran}}\\
 {\small{ $^{{\rm b}}$School of Mathematics, Institute for Research in Fundamental Sciences
 (IPM)}}\vspace{-1mm}\\ {\small{ P.O.Box: 19395-5746, Tehran, Iran}}\vspace{-1mm}\\
 {\small{mbehbood@cc.iut.ac.ir}}\vspace{-1mm}\\
 {\small{hshojaee@math.iut.ac.ir}}}
\date{}
\begin{document}
\maketitle
\begin{abstract}
\noindent{ A well-known result of   K\"{o}the and Cohen-Kaplansky
states that a commutative ring $R$ has the property that every
$R$-module is a direct sum of cyclic modules if and only if $R$ is
an Artinian principal ideal ring.
 This motivated us to study commutative rings for which every ideal is a direct sum of cyclic
modules. Recently, in  [M. Behboodi, A. Ghorbani, A.
Moradzadeh-Dehkordi, {\it Commutative Noetherian local rings whose
ideals are  direct sums  of cyclic modules,} J. Algebra  345
(2011) 257--265] the authors considered this question in the
context of finite direct products of commutative Noetherian local
rings. In this paper, we continue their study by dropping the
Noetherian condition. }
\end{abstract}

{\bf ~~~~~~~~~~~~~~~~~~~~~~~~~~~~~~~~1. Introduction}

\noindent The study of rings over which modules  are direct sums
of cyclic modules has a long history.
 The first
important contribution in this direction is due to K\"{o}the
\cite{Kothe} who considered rings over which all  modules are
direct sums of cyclic modules. K\"{o}the showed that over an
Artinian principal ideal ring, each module is a direct sum of
cyclic modules. Furthermore, if a commutative Artinian ring has
the property that all its modules are direct sums of cyclic
modules, then it is necessarily a principal ideal ring. Later,
Cohen and Kaplansky \cite{Cohen} obtained  the following.

 \noindent {\bf Result 1.1.} (Cohen and Kaplansky, \cite{Cohen})
{\it If $R$ is a commutative ring such that each $R$-module is a
direct sum of cyclic modules, then $R$ must be an Artinian
principal ideal ring.}

An interesting natural question arises. Instead of considering
rings for which {\it all } modules are direct sums of cyclic
modules, we weaken this condition and study rings $R$ for which it
is assumed only that the {\it ideals} of $R$ are direct sums of
cyclic modules. The study of such commutative rings was initiated
by  Behboodi, Ghorbani and Moradzadeh-Dehkordi in
\cite{Behboodi1}. In particular, they established the following
theorem.

\noindent {\bf Result 1.2.} (\cite[Theorem 2.11]{Behboodi1}) {\it
Let $(R, \cal{M})$  be a commutative Noetherian local  ring, where
${\mathcal{M}}$ denotes the unique maximal  ideal of $R$. Then the
following
statements are equivalent:}\vspace{2mm}\\
(1) {\it Every ideal of $R$ is a direct sum of cyclic $R$-modules.}\vspace{1mm}\\
(2) {\it There exist  an positive integer  $n$ and a set of
elements $\{w_1,\cdots, w_n\}\subseteq  R$ such that
 \indent  ${\cal{M}}=Rw_1\oplus\cdots\oplus Rw_n$ with at most two
  of $R{w_i}^,s$  not  simple.}\vspace{1mm}\\
(3) {\it There exists  an positive integer  $n$ such that every
ideal of $R$ is a direct sum of at most \indent $n$
 cyclic   $R$-modules.}\vspace{1mm}\\
(4) {\it Every ideal of $R$ is a direct summand of a direct sum of
cyclic $R$-modules.}

In this paper  we consider commutative local rings for which every
ideal is a direct sum of cyclic modules, that is, we drop the
Noetherian condition from [1]. In particular, we describe the
ideal structure of such rings.

In the sequel  all rings are commutative with identity and all
modules are unital. For a ring $R$, we denote (as usual) the set
of prime ideals of $R$ by Spec$(R)$. Also, Nil$(R)$ is the ideal
of all nilpotent elements of $R$. We denote the classical Krull
dimension of $R$ by dim$(R)$. Let $X$ be either an element or a
subset of $R$. The  annihilator of $X$ is the ideal Ann$(X)=\{a\in
R~|~ aX=0\}$. A ring $R$  is {\it local} in case $R$ has a unique
maximal ideal. In this paper  $(R,{\mathcal{M}})$ will be a local
ring with maximal ideal  ${\mathcal{M}}$. An $R$-module $N$ is
called {\it simple} if $N\neq (0)$ and it has no submodules except
$(0)$ and $N$. An $R$-module $M$ is a {\it semisimple}  module if
it is a direct sum of simple modules. Also, an $R$-module $M$ is
called a {\it homogenous semisimple} $R$-module if it is a direct
sum of isomorphic simple $R$-modules, i.e., ${\rm Ann}(M)$ is a
maximal ideal of $R$.

 It will be shown (see Theorems 3.1
and 3.3) that if a local ring $(R, \cal{M})$ has the property that
every   ideal of $R$  is a direct sum of cyclic $R$-modules, then
dim$(R)\leq 1$ and $|{\rm{Spec}}(R)|\leq 3$. Moreover, there is an
index set $\Lambda$ and a set of elements $\{x, y\} \cup
\{w_{\lambda}\}_{\lambda\in \Lambda}\subseteq R$ such that
${\cal{M}}=Rx\oplus Ry\oplus(\bigoplus_{\lambda\in \Lambda}
Rw_{\lambda})$ with: each $Rw_{\lambda}$ a simple $R$-module,
$R/{\rm Ann}(x)$, $R/{\rm Ann}(y)$  principal ideal rings, and
${\rm{Spec}}(R)\subseteq \{(0), {\cal{M}}, Rx \oplus
(\bigoplus_{\lambda\in \Lambda}Rw_{\lambda}),Ry \oplus
(\bigoplus_{\lambda\in \Lambda}Rw_{\lambda})\}$. Also, we prove
the following main theorem.

\noindent {\bf Result 1.3.} (See Theorem 3.7) {\it For a local
ring $(R, \cal{M})$ the following
statements are equivalent:}\vspace{2mm}\\
(1)  {\it Every ideal of $R$ is a direct sum of cyclic $R$-modules.}\vspace{1mm}\\
(2)  {\it Every ideal of $R$ is  a direct sum of cyclic
$R$-modules, at most two of which
are not \indent simple.}\vspace{1mm}\\
(3) {\it  There is an index set $\Lambda$ and a set of elements
$\{x, y\} \cup \{w_{\lambda}\}_{\lambda\in \Lambda}\subseteq R$
such that ${\cal{M}}= \indent Rx  \oplus  Ry
\oplus(\bigoplus_{\lambda\in \Lambda} Rw_{\lambda})$ with: each
$Rw_{\lambda}$ a simple $R$-module, $R/{\rm Ann}(x)$ and $R/{\rm
Ann}(y)$ \indent principal ideal
rings.}\vspace{1mm}\\
(4) {\it  Every ideal of $R$ is a direct summand of a direct sum
of cyclic $R$-modules.}

Finally, some relevant
examples and counterexamples are indicated in Section 4.\\

\noindent{\bf 2. Preliminaries}

\noindent We begin  this section with the following result from
commutative algebra due to I. M. Isaacs which states that to check
whether every ideal in a ring  is principal, it suffices to test only
the prime ideals.

\noindent{\bf Lemma 2.1.} (Attributed to I. M. Isaacs in  \cite[p.
8, Exercise 10]{Kaplansky}) {\it A commutative ring $R$  is a
principal ideal ring if and only if every prime ideal of $R$ is a
principal ideal.}

\noindent{\bf Proposition 2.2.} {\it Let $R$  be a  ring. If every
prime ideal of $R$ is a direct sum of cyclic $R$-modules, then
$R/P$ is a principal ideal domain {\rm (PID)} for each prime ideal
$P\in {\rm Spec}(R)$. Consequently, ${\rm dim}(R)\leq 1$.}

\noindent{\bf Proof.} Assume that $P\in {\rm Spec}(R)$ and $Q/P$
is a prime ideal of $R/P$. Since $Q\in {\rm Spec}(R)$,
$Q=\oplus_{i\in I} Rx_i $ for some index set $I$ and $x_i\in R$
for each $i\in I$. If $Q/P$ is nonzero, then there exists $j\in I$
such that $x_j\notin P$. Since for each $i\in I$,
$Rx_iRx_j=(0)\subseteq P$, we conclude that $Rx_i\subseteq P$ for
each $i\neq j$. It follows that $x_j+P\in R/P$ is a generator  for
$Q/P$. Thus  by Lemma 2.1, $R/P$ is a PID. Since this holds for
all prime ideals $P$ of $R$, we conclude that dim$(R)\leq
1$.$~\square$

The following two results from \cite{Behboodi1}  are crucial to
our investigation.

\noindent{\bf Lemma 2.3.} (\cite[Proposition 2.2]{Behboodi1}) {\it
Let $(R, \cal{M})$ be a local ring. Suppose that
${\cal{M}}=Rx\oplus Ry \oplus Rz\oplus K$ for nonzero elements
$x$, $y$, and $z$ and an ideal $K$ of $R$. Further suppose that
neither  of  $Rx$, $Ry$, or $Rz$ is a simple $R$-module. Then the
ideal $J:=R(x+y)+R(x+z)$ is not a direct sum of cyclic
$R$-modules.}

\noindent{\bf Lemma 2.4.} (\cite[Corollary  2.3]{Behboodi1}) {\it
Suppose $(R, \cal{M})$ is a local ring such that every ideal of
$R$ is a direct sum of cyclic $R$-modules. Then there is an index
set $\Lambda$ and a set of elements $\{w_{\lambda}\}_{\lambda\in
\Lambda}\subseteq R$ such that ${\cal{M}}= \bigoplus_{\lambda\in
\Lambda} Rw_{\lambda}$ with at most two of the ${Rw_{\lambda}}$'s
not simple.}

 \noindent{\bf Lemma 2.5.} (See  \cite[Proposition 3]{Warfield1}) {\it Let $R$ be a local ring and
 let $M$ be an  $R$-module. If there is an index set
$\Lambda$ and a set of ideals $\{I_{\lambda}\}_{\lambda\in
\Lambda}$ such that  $M=\bigoplus_{\lambda\in \Lambda}
R/I_{\lambda}$, then every direct summand of $M$ is also a direct
sum of cyclic $R$-modules, each isomorphic to one of the
$R/I_{\lambda}$.}

\noindent{\bf Lemma 2.6.} {\it  Let $R$ be a ring and let $M$ be a
homogenous semisimple $R$-module.  Then  there is an index set
$\Lambda$ and a set of elements $\{w_{\lambda}\}_{\lambda\in
\Lambda}\subseteq R$ such that $M= \bigoplus_{\lambda\in \Lambda}
Rw_{\lambda}$ where $Rw_{\lambda}$'s are isomorphic simple
$R$-modules
 and  every  submodule of $M$ is  also of the form
$N=\bigoplus_{\gamma\in\Gamma}Rw'_{\gamma}$ where $\Gamma$ is an
index set with $|\Gamma|\leq|\Lambda|$ and $Rw'_{\gamma}$'s are
isomorphic simple $R$-modules.}

\noindent{\bf Proof.} The proof is clear from the fact that ${\rm
Ann}(M)$ is a maximal ideal of $R$ and
   $M$ is an $R/{\rm Ann}(M)$-vector space.$~\square$

We conclude this section with the following proposition from
\cite{Behboodi1} that provides an analogue of the Invariant Base
Number of a free module over a commutative ring.

 \noindent{\bf Lemma 2.7.} (\cite[Proposition 2.15]{Behboodi1}) {\it Let $R$ be a ring.
 The following statements are equivalent:}\vspace{2mm}\\
 (1) {\it $R$ is a local ring.}\vspace{1mm}\\
  (2) {\it  If $\bigoplus_{i=1}^nRx_i\cong \bigoplus_{j=1}^mRy_j$ where $n$, $m\in \Bbb{N}$
   and $\forall i,j$,
 $Rx_i$,  $Ry_j$  are nonzero cyclic  \indent $R$-modules,   then  $n=m$.}\vspace{1mm}\\
 (3) {\it  If $\bigoplus_{i\in I}Rx_i\cong \bigoplus_{j\in J}Ry_j$ where $I$, $J$ are index sets and
 $Rx_i$,  $Ry_j$  are nonzero cyclic \indent  $R$-modules, then
 $|I|=|J|$.}\\

\noindent{\bf 3. Main Results}\vspace{3mm}

\noindent{\bf Theorem  3.1.} {\it Let $(R, \cal{M})$ be a
 local ring such that every ideal of $R$ is a direct sum of
cyclic $R$-modules. Then {\rm dim}$(R)\leq 1$,  $|{\rm
Spec}(R)|\leq 3$, and there is an index set $\Lambda$ and a set of
elements $\{x, y\} \cup \{w_{\lambda}\}_{\lambda\in
\Lambda}\subseteq R$ such that ${\cal{M}}=Rx  \oplus  Ry
\oplus(\bigoplus_{\lambda\in \Lambda} Rw_{\lambda})$ with: each
$Rw_{\lambda}$ a simple $R$-module.
  Moreover,} \vspace{1mm}\\
(a)  {\it If $x$, $y\in {\rm Nil}(R)$, then
${\rm Spec}(R)=\{{\cal{M}}\}$.}\vspace{1mm}\\
(b)  {\it If ${\cal{M}}=Rz$ and $z\not\in {\rm Nil}(R)$, then
 ${\rm Spec}(R)=\{(0), {\cal{M}}\}$.}\vspace{1mm}\\
(c)  {\it If ${\cal{M}}$ is not cyclic,  $x\not\in {\rm Nil}(R)$
and $y\in {\rm Nil}(R)$, then ${\rm Spec}(R)=\{{\cal{M}},
Ry\oplus(\bigoplus_{\lambda\in \Lambda} Rw_\lambda)\}$.}\vspace{1mm}\\
(d) {\it If ${\cal{M}}$ is not cyclic,  $x\in {\rm Nil}(R)$ and
$y\not\in {\rm Nil}(R)$, then ${\rm Spec}(R)=\{{\cal{M}},
Rx\oplus(\bigoplus_{\lambda\in \Lambda} Rw_\lambda)\}$.}\vspace{1mm}\\
(e) {\it If ${\cal{M}}$ is not cyclic and   $x,~y\not\in
{\rm Nil}(R)$, then\vspace{1mm}\\
 $~~~~~~~~~~~~~~~~~~~~~~~~~{\rm Spec}(R)=\{{\cal{M}},
Rx\oplus(\bigoplus_{\lambda\in \Lambda} Rw_\lambda),
Ry\oplus(\bigoplus_{\lambda\in \Lambda} Rw_\lambda)\}$. }

\noindent{\bf Proof.}  By Proposition 2.2, dim$(R)\leq 1$. Also,
by
 Lemma 2.4, there is an index set $\Lambda$ and a set of elements $\{x, y\}
\cup \{w_{\lambda}\}_{\lambda\in \Lambda}\subseteq R$ such that
${\cal{M}}=Rx  \oplus  Ry \oplus(\bigoplus_{\lambda\in \Lambda}
Rw_{\lambda})$ with: each $Rw_{\lambda}$ a simple $R$-module. Thus
by \cite[Lemma 2.1]{Behboodi1}, $w^2_{\lambda}=0$ for each
${\lambda\in \Lambda}$. We consider the following five cases.

   Case (a): Suppose that $x$, $y\in{\rm Nil}(R)$. Since
$w_{\lambda}^2=0$ for each ${\lambda}\in {\Lambda}$, we conclude
that Nil$(R)={\cal{M}}$ and thus  ${\rm{Spec}}(R)=\{{\cal{M}}\}$.

  Case (b): Suppose that ${\cal{M}}=Rz$ and $z\not\in{\rm Nil}(R)$.
Then  ${\rm {dim}(R)}=1$. Let $P\in {\rm{Spec}(R)}\setminus
\{{\cal{M}}\}$. Since $P\subsetneq {\cal{M}}=Rz$,   $P=Pz$ and,
 by  Nakayama's lemma,  $P=(0)$  (since Nakayama's lemma holds for
any direct sum of finitely generated modules).
  Thus, by Lemma 2.1, $R$ is a principal ideal domain and
${\rm{Spec}}(R)=\{{(0), \cal{M}}\}$.

  Case (c): Suppose that ${\cal{M}}$ is not cyclic, $x\not\in{\rm
Nil}(R)$ and $y\in{\rm Nil}(R)$. Then Nil$(R)\neq {\cal{M}}$ and
dim$(R)=1$. Let $P$ be a prime ideal of $R$ such that $P\subsetneq
{\cal{M}}$. Since $w^2_\lambda=0$ for each ${\lambda}\in
{\Lambda}$, we conclude that $Ry\oplus (\bigoplus_{{\lambda}\in
{\Lambda}}Rw_{\lambda})\subseteq P$. Thus $P=(Rx\cap P)\oplus
Ry\oplus (\bigoplus_{{\lambda}\in {\Lambda}}Rw_{\lambda})$. Since
${\cal{M}}\not\subseteq P$, it follows that $x\not\in P$ and so
$Rx\cap P=Px$. Thus $P=Px\oplus Ry\oplus(\bigoplus_{{\lambda}\in
{\Lambda}}Rw_{\lambda})$ and hence, $Px=Px^2=RxPx$. By Nakayama's
lemma,  $Px=0$. Thus $P=Ry\oplus (\bigoplus_{{\lambda}\in
{\Lambda}}Rw_{\lambda})$.  Therefore, ${\rm{Spec}}(R)=\{{\cal{M}},
Ry\oplus (\bigoplus_{{\lambda}\in {\Lambda}}Rw_{\lambda})\}$.

  Case (d): If ${\cal{M}}$ is not cyclic, $x\in{\rm Nil}(R)$ and
$y\not\in{\rm Nil}(R)$ and a similar argument allows us to
conclude that ${\rm{Spec}}(R)=\{{\cal{M}}, Rx\oplus
(\bigoplus_{{\lambda}\in {\Lambda}}Rw_{\lambda})\}$.

     Case (e): Suppose that ${\cal{M}}$ is not cyclic and  $x$, $y\not\in {\rm
Nil}(R)$. Thus Nil$(R)\neq {\cal{M}}$ and so dim$(R)=1$. Let $P$
be a prime ideal of $R$ such that $P\subsetneq {\cal{M}}$.  Since
$xy=0$, $x\in P$ or $y\in P$. If $x\in P$, then $y\not\in P$ and
$P=Rx\oplus(P\cap Ry)\oplus(\bigoplus_{{\lambda}\in
{\Lambda}}Rw_{\lambda})$. As in Case (c), we have  $Ry\cap
P=Py=(0)$ and
 $P=Rx\oplus
(\bigoplus_{{\lambda}\in {\Lambda}}Rw_{\lambda})$. Similarly, if
$y\in P$, then $P=Ry\oplus (\bigoplus_{{\lambda}\in
{\Lambda}}Rw_{\lambda})$. On the other hand, since $x$, $y\not\in
{\rm Nil}(R)$, there exist $P_1$, $P_2\in {\rm{Spec}(R)}\setminus
\{\cal{M}\}$ such that $x\in P_1$, $y\not\in P_1$ and $x\not\in
P_2$, $y\in P_2$. Therefore, ${\rm{Spec}}(R)=\{{\cal{M}}, Rx
\oplus (\bigoplus_{{\lambda}\in {\Lambda}}Rw_{\lambda}),Ry \oplus
(\bigoplus_{{\lambda}\in {\Lambda}}Rw_{\lambda})\}$.$~\square$

 We can now state the following corollary,
 an analog of both Kaplansky's theorem \cite[Theorem
12.3]{Kaplansky1}  (which states that a commutative Noetherian
ring $R$ is a principal ideal ring if and only if every maximal
ideal of $R$ is principal) and  Cohen's theorem [2] (which states
that $R$ is Noetherian if and only if every prime ideal of $R$ is
finitely generated).

 \noindent{\bf Corollary  3.2.} {\it Let $(R, \cal{M})$ be a
 local ring such that every ideal of $R$ is a direct sum of
cyclic $R$-modules. Then $\cal{M}$ is  cyclic (resp. finitely
generated) if and only if  $R$ is a principal ideal ring (resp.
Noetherian ring).}

\noindent{\bf Proof.} From Theorem 3.1, we see that if $\cal{M}$
is principal (resp. finitely generated), then the same holds for
every prime ideal of $R$. The conclusion follows from Lemma 2.1
(resp. Cohen's theorem).$~\square$

Next, we sharpen Corollary 2.3 of \cite{Behboodi1} (see Lemma
2.4).

\noindent{\bf Theorem  3.3.} {\it Let $(R, \cal{M})$  be a local
ring such that every ideal of $R$ is a direct sum of cyclic
 $R$-modules. Then, with the notation
of Theorem 3.1, both  $R/{\rm Ann}(x)$ and  $R/{\rm Ann}(y)$ are
principal ideal rings. }

\noindent{\bf Proof.} Using Theorem 3.1, we see that each prime
ideal of $S:=R/(Ry\oplus(\bigoplus_{\lambda\in \Lambda}
Rw_{\lambda}))$ is principal. Now Lemma 2.1 implies that $S$ is a
principal ideal ring. Since $Ry\oplus (\bigoplus_{\lambda\in
\Lambda} Rw_{\lambda})\subseteq  {\rm Ann}(x)$, $R/ {\rm Ann}(x)$
is a homomorphic image of $S$ and  is thus a principal ideal ring.
Similarly, $R/ {\rm Ann}(y)$ is a principal ideal ring.$~\square$

Let us now outline the proof of the main theorem of this paper.
We have divided the proof into a sequence of
propositions.

\noindent{\bf Proposition 3.4.} {\it Let  $(R, \cal{M})$ be a
local ring such that  ${\cal{M}} = Rx \oplus L$ for some ideal $L$
of $R$ and some $x\in  R$. Suppose further that $R/{\rm Ann}(x)$
is a principal ideal ring. Then every nonzero element of $Rx$ is
of the form $ax^n$ for some unit $a$ and positive integer $n$.}

\noindent{\bf Proof.} We identify  $Rx$ with $\bar{R}:=R/{\rm
{Ann}}(x)$. Then  an element $a\in R$ is a unit whenever $\bar{a}$
is  a unit in ${\bar{R}}$ (in fact,  if $\bar{a}\bar{b}=\bar{1}$
for some $b\in R$, then $ab-1\in {\rm {Ann}}(x)\subseteq{\cal{M}}$
and so $a$ is a unit in $R$). Since ${\bar{R}}$ is a local
principal ideal ring, with maximal ideal ${\bar{R}}{\bar{x}}$,
${\bar{R}}$ is a Noetherian ring and so  by Krull's intersection
theorem $\cap_{i=1}^{\infty}{\bar{R}}{\bar {x}}^i=(0)$. Suppose
that $0\neq rx\in Rx$ where $r$ is not a unit. Thus $0\neq
\bar{r}\in {\bar{R}}{\bar{x}}$. We claim that $\bar{r}={\bar
{a}}{\bar{x}}^n$ for some unit ${\bar {a}}$ and $n\in\Bbb{N}$. If
not,   for each $i\in\Bbb{N}$ there exists ${\bar{r}}_i\in
{\bar{R}}{\bar {x}}$ such that ${\bar{r}}={\bar{r}}_i{\bar{x}}^i$,
i.e.,  ${\bar{r}}\in \cap_{i=1}^{\infty}{\bar{R}}{\bar
{x}}^i=(0)$, a contradiction. Therefore, $\bar{r}={\bar
{a}}{\bar{x}}^n$ for some unit ${\bar {a}}$ and $n\in\Bbb{N}$.
Consequently, $a$ is a unit in $R$ and $rx=ax^{n+1}$,   which
completes the proof.$~\square$

 \noindent{\bf Proposition 3.5.} {\it Let $(R, \cal{M})$ be a
local ring such that  there is an index set $\Lambda$ and a set of
elements $\{x\} \cup \{w_{\lambda}\}_{\lambda\in \Lambda}\subseteq
R$ such that ${\cal{M}}=Rx \oplus(\bigoplus_{\lambda\in \Lambda}
Rw_{\lambda})$ with: each $Rw_{\lambda}$ a simple $R$-module and
$R/ {\rm Ann}(x)$ principal ideal ring.
 Then every proper  ideal of
$R$ is of the form $I=Rx'\oplus(\bigoplus_{\gamma\in \Gamma}
Rw'_{\gamma})$ where $\Gamma$ is an index set with
$|\Gamma|\leq|\Lambda|$, ${Rw'_{\gamma}}$'s ($\gamma\in \Gamma$)
are simple $R$-modules   and $x'\in R$ is such that $R/{\rm
Ann}(x')$ is a principal ideal ring.}

\noindent{\bf Proof.}    Assume that $I$ is a proper ideal of $R$
and  $L=\bigoplus_{\lambda\in \Lambda} Rw_{\lambda}$. Clearly, $L$
is a homogenous semisimple  $R$-module with $L^2=(0)$. If
$I\subseteq Rx$, then every ideal contained in $I$ is principal,
and we are done since $R/{\rm Ann}(x)$ is a principal ideal ring.
If $I\subseteq L$, then by Lemma 2.6,  $I$ is a direct sum of at
most $|\Lambda|$ simple  modules. Thus we can assume that
$I\not\subseteq Rx$, $I\not\subseteq L$ and $(0)\subsetneq
I\subsetneq {\cal{M}}$. By Proposition 3.4, there exist $n\in
\Bbb{N}$ and $l\in L$ such that $x^n+l\in I$. Among all such
expressions, choose one, $x^{n_0}+l_{0}$, with $n_0$ minimal. We
set $x'=x^{n_0}$ if $x^{n_0}\in I$, otherwise we set
$x'=x^{n_0}+l_0$. Set
 $$J:=\{l\in L~| ~ax+l\in I,~ {\rm for~ some~} a\in R\}.$$ Then
$J$ is an ideal of $R$ with $J\subseteq L$.  We next prove that
 $Rx'\cap (I\cap J)=(0)$. For see this, let  $rx'\in Rx'\cap (I\cap J)$ where  $r\in R$. If
$x^{n_0}\in I$, then  $x'=x^{n_0}$ and so $rx'=rx^{n_0}\in Rx\cap
 L=(0)$. If   $x^{n_0}\not\in I$, then  $rx'=r(x^{n_0}+l_0)=l$ for some $l\in
I\cap J\subseteq L$. Hence we have  $rx^{n_0}=l-rl_0\in Rx\cap
L=\{0\}$. Since $x^{n_0}\not\in I$,  $r$ is not a unit and so by
Proposition 3.4, $r=ax^n+l_1$ where $a$ is a unit, $l_1\in L$ and
$n\in\Bbb{N}$. Since $L^2=(0)$, we conclude that
$rx'=(ax^n+l_1)(x^{n_0}+l_0)=ax^{n+n_0}=l\in Rx\cap L=\{0\}$. Thus
 $Rx'\cap (I\cap J)=(0)$ and $Rx'\oplus (I\cap J)\subseteq I$. In
fact, we will show that $I= Rx'\oplus (I\cap J)$. Assume that
$u=ax^s+l\in I$ where $a\in (R\setminus {\cal{M}})\cup \{0\}$,
$s\in \Bbb{N}$ and $l\in L$. If $a=0$, then $l\in I\cap J$ and so
$u=l\in Rx'\oplus (I\cap J)$. Thus we can assume that $a\in
R\setminus {\cal{M}}$. Therefore,
$$u-ax^{s-n_{0}}x'=(ax^{s}+l)-ax^{s-n_{0}}(x^{n_0}+l_{0})=
 l+ax^{s-n_{0}}l_{0}\in I.$$  Since $l$ and $l_0$ contained in  $J$,
  $l+ax^{s-n_{0}}l_{0}\in J$. Hence it follows that $u\in Rx'\oplus
  (I\cap J)$.
  Therefore, $I= Rx'\oplus (I\cap J)$.  By
Lemma 2.6,  $I\cap L$ is a direct sum of at most $|\Lambda|$
simple $R$-modules. Since ${\cal{M}}l_0=(0)$, we conclude that
${\rm Ann}(x^{n_0}+l_0)= {\rm Ann}(x^{n_0})$. This implies that
$Rx'\cong Rx^{n_0}\subseteq Rx$. Since $R/{\rm Ann}(x)$ is a
principal ideal ring,
 $R/ {\rm Ann}(x')$ is also a principal ideal ring.$~\square$

 \noindent{\bf Proposition  3.6.} {\it Let
$(R, \cal{M})$ be a local ring such that there is an index set
$\Lambda$ and a set of elements $\{x, y\} \cup
\{w_{\lambda}\}_{\lambda\in \Lambda}\subseteq R$ such that
${\cal{M}}=Rx  \oplus  Ry \oplus(\bigoplus_{\lambda\in \Lambda}
Rw_{\lambda})$ with: each $Rw_{\lambda}$ a simple $R$-module and
$R/ {\rm Ann}(x)$,  $R/ {\rm Ann}(y)$  principal ideal rings. If
$I$ is an ideal of $R$, then one of the following
holds:}\vspace{2mm}\\
(i) {\it There is an index set $\Gamma$ and a set of elements
$\{x', y'\} \cup \{w'_{\gamma}\}_{\gamma\in \Gamma}\subseteq R$
such that $I= \indent Rx'  \oplus  Ry' \oplus(\bigoplus_{\gamma\in
\Gamma} Rw'_{\gamma})$ with: $|\Gamma|\leq|\Lambda|$, each
$Rw'_{\gamma}$ a simple $R$-module, and  $R/ {\rm Ann}(x')$,
\indent
 $R/ {\rm Ann}(y')$
principal ideal rings.}\\
(ii) {\it There is an index set $\Gamma$ and a set of elements
$\{z'\} \cup \{w'_{\gamma}\}_{\gamma\in \Gamma}\subseteq R$ such
that $I= \indent Rz' \oplus(\bigoplus_{\gamma\in \Gamma}
Rw'_{\gamma})$ with: $|\Gamma|\leq|\Lambda|$, each $Rw'_{\gamma}$
a simple $R$-module, and each ideal of  \indent $R/{\rm Ann}(z')$
a direct sum of at most  two principal ideals.}

\noindent{\bf Proof.}  By Proposition 3.5  we can assume that $Rx$
and $Ry$ are not simple $R$-modules. Let $L=\bigoplus_{\lambda\in
 \Lambda}Rw_{\lambda}$. Clearly both $R/Rx$ and $R/Ry$ are local rings with
maximal ideals  ${\cal{M}}_{x}={\cal{M}}/Rx\cong Ry\oplus  L$ and
${\cal{M}}_{y}={\cal{M}}/Ry\cong Rx\oplus  L$.
 Assume that $I$
is an ideal of $R$. First, note that if $I\subseteq Rx\oplus L$
(resp. $I\subseteq Ry\oplus L$), then $I\cong (I\oplus Ry)/Ry$
(resp. $I\cong (I\oplus Rx)/Rx$) and  so, by Proposition 3.5,
$I=Rx'\oplus(\bigoplus_{\gamma\in \Gamma} Rw'_{\gamma})$ where
$\Gamma$ is an index set with $|\Gamma|\leq|\Lambda|$, the
${Rw'_{\gamma}}$'s ($\gamma\in \Gamma$) are simple $R$-modules,
and $x'\in R$ is such that $R/{\rm Ann}(x')$ is a principal ideal
ring. Therefore, according to  the above remark, we  can assume
that $I\not\subseteq Rx\oplus L$, $I\not\subseteq Ry\oplus L$ and
$(0)\subsetneq I \subsetneq {\cal{M}}$.

By Proposition  3.4, every element of $I$ has the form
$ax^s+by^t+l$ where $a\in (R\setminus {\cal{M}})\cup \{0\}$, $b\in
(R\setminus {\cal{M}})\cup\{0\}$, $s,~t\in \Bbb{N}$ and $l\in L$.
Since $I\not\subseteq Rx\oplus L$, and $I\not\subseteq Ry\oplus
L$, it follows that there exist $e_1$, $e_2\in I$ where
 $e_1=a_{1}x^n+b_{1}y^t+l_1$, and
$e_2=a_{2}x^s+b_{2}y^m+l_2$  for some   $a_1, ~b_2\in R\setminus
{\cal{M}}, ~a_2,~b_1\in R $,
 $l_1,~l_2\in L$ and $n,~s,~t,~m\in \Bbb{N}$ (in fact $e_1\in I\setminus Rx\oplus L$ and
 $e_2\in I\setminus Ry\oplus L$).
  Thus $xe_1=a_1x^{n+1}\in I$ and $ye_2=b_2y^{m+1}\in I$.
Since  $a_1, ~b_2\in R\setminus {\cal{M}}$,  $x^{n+1}\in I$ and
$y^{m+1}\in I$. Suppose that $n_0$ (resp. $m_0$) is  the smallest
natural number  such that  $x^{n_0}+l_1\in I$ (resp.
$y^{m_0}+l_2\in I$) for some $l_1\in L$ (resp. $l_2\in L$). We set
$x'=x^{n_0}$ if  $x^{n_0}\in I$, otherwise we set
$x'=x^{n_0}+l_1$. Also, we set $y'=y^{m_0}$ if $y^{m_0}\in I$,
otherwise we set $y'=y^{m_0}+l_2$.   Set
 $$J:=\{l\in L~| ~ax+by+l\in I,~ {\rm for~ some~} a,~b\in R\}.$$ Then
$J$ is an ideal of $R$ with $J\subseteq  L$. On can easily see
that  $Rx'+ Ry'+ (I\cap J)=Rx'\oplus Ry'\oplus (I\cap J)\subseteq
I$.
 Now we proceed by cases.

  {Case  (a)}: Suppose that $I=Rx'\oplus Ry'\oplus (I\cap J)$. By Lemma 2.6, $I\cap J$ is a direct sum of
   at most $|\Lambda|$ simple $R$-modules.  It follows that $I$ is
   a direct sum of at most $|\Lambda|+2$ cyclic modules by Lemma 2.7. Also, if $x'\neq 0$, then  $x^{n_0}\neq 0$ and so
${\rm Ann}(x')={\rm Ann}(x^{n_0})$ since ${\cal{M}}l_1=(0)$. It
follows that $R/{\rm Ann}(x')$ is a principal ideal ring (since
$Rx^{n_0}\subseteq Rx$). We conclude similarly that either
$Ry'=(0)$ or $R/{\rm Ann}(y')$ is a principal ideal ring.

 {Case  (b)}: Suppose that $Rx'\oplus Ry'\oplus (I\cap J)\subsetneq I$. We claim that
 every element of $I\setminus (Rx'\oplus Ry'\oplus (I\cap J))$
 is of the form $cx^{n_0-1}+dy^{m_0-1}+l$
 where $c$, $d\in R\setminus {\cal{M}}$ and $l\in L$. Let
  $z=cx^{s}+dy^t+l\in I\setminus (Rx'\oplus Ry'\oplus (I\cap
J))$ where $c,d\in (R\setminus {\cal{M}})\cup \{0\}$ and $l\in L$.
If $c=0$, then $z=dy^t+l\in I$ and so $t\geq m_0$. If $t> m_0$,
then $z=dy^{t-m_0}(y^{m_0}+l_2)+l\in Ry'\oplus (I\cap J)$, a
contradiction. Thus $t=m_0$ and this implies that
$z=dy'+(l-dl_1)\in Ry'\oplus (I\cap J)$, a contradiction. Thus
$c\in R\setminus {\cal{M}}$.  We conclude similarly that $d\in
R\setminus {\cal{M}}$. Our next  claim  is that $s=n_0-1$ and
$t=m_0-1$. If not, to obtain a contradiction, as we see by
considering the following three subcases:

{Subcase  (i)}: Suppose that $s< n_0-1$ or $t< m_0-1$. There is no
loss of generality in assuming that $s< n_0-1$. Then
$x^{n_0-1}=x^{n_0-1-s}c^{-1}z\in I$  which contradicts the
minimality of $n_0$.

{Subcase (ii)}: Suppose that $s\geq n_0$ and $t\geq m_0$. Then
$z=cx^{s-n_0}x'+dy^{t-m_0}y'-(cx^{s-n_0}l_1-dy^{t-m_0}l_2)+l\in
I$. Since $z-cx^{s-n_0}x'+dy^{t-m_0}y'\in I$,
$(cx^{s-n_0}l_1-dy^{t-m_0}l_2)+l\in (I\cap J)$ and hence $z\in
Rx'\oplus Ry'\oplus (I\cap J) $, a contradiction.

{Subcase (iii)}: Suppose that $s\geq n_0$,  $t= m_0-1$ or
$s=n_0-1$, $t\geq m_0$. Without loss  of generality we can assume
$s\geq n_0$ and $t= m_0-1$.  Then
$z=cx^{s-n_{0}}x'+dy^{m_0-1}-cx^{s-n_{0}}l_1+l$ and so
$z-cx^{s-n_{0}}x'=dy^{m_0-1}-cx^{s-n_{0}}l_1+l\in I$. Thus
$y^{m_0-1}+d^{-1}(l-cx^{s-n_{0}}l_1)\in I$ which contradicts the
minimality of $m_0$.

Therefore, every element of $I\setminus (Rx'\oplus Ry'\oplus
(I\cap J))$ has the form $cx^{n_0-1}+dy^{m_0-1}+l$ where $c,d\in
R\setminus {\cal{M}}$ and $l\in L$. Let
$z'=cx^{n_0-1}+dy^{m_0-1}+l\in I\setminus (Rx'\oplus Ry'\oplus
(I\cap J))$ where $c,d\in R\setminus {\cal{M}}$ and $l\in L$.
Since $L^2=(0)$, it is easy to check that $Rz'\cap (I\cap J)=(0)$.
We now claim that $I=Rz'\oplus (I\cap J)$. Note that if
$x'=x^{n_0}$, then $x'=xc^{-1}z'$ and so $x'\in Rz'$. Also, if
$x'=x^{n_0}+l_1$ for some nonzero element $l_1$ of $L$, then
$x'=xc^{-1}z'+l_1$. Therefore, since $l_1\in J$ and
$l_1=x'-xc^{-1}z'\in I$, so $l_1\in I\cap J$. Hence $x'\in
Rz'\oplus (I\cap J)$. We conclude similarly that $y'\in Rz'\oplus
(I\cap J)$. Thus $Rx'\oplus Ry'\oplus (I\cap J)\subseteq Rz'\oplus
(I\cap J)\subseteq I$. Suppose, contrary to our claim, that
$I\not\subseteq Rz'\oplus (I\cap J)$. Then there exists an element
$u\in I\setminus (Rz'\oplus (I\cap J))$ and so $u\in I\setminus
(Rx'\oplus Ry'\oplus (I\cap J))$. Therefore,
$u=c'x^{n_0-1}+d'y^{m_0-1}+l'$ for some  $c'$, $d'\in R\setminus
{\cal{M}}$ and $l'\in L$. Then
$$u-c'c^{-1}z'=(c'x^{n_0-1}+d'y^{m_0-l}+l')-c'c^{-1}(cx^{n_0-1}+dy^{m_0-1}+l)=d''y^{m_0-1}+l''\in I$$
where $d''=d'-c'c^{-1}d$ and $l''=l'-c'c^{-1}l$. If $d''=0$, then
$l''\in (I\cap J)$ and so $u=c'c^{-1}z'+l''\in Rz'\oplus (I\cap
J)$, a contradiction.
  If $d''\in R\setminus {\cal{M}}$, then $y^{m_0-1}+d''^{-1}l''\in I$ which contradicts
 the minimality of $m_0$.  Thus  $d''\in
{\cal{M}}$ and  there exists $r\in R$ such that
$d''y^{m_0-1}=ry^{m_0}$.  Therefore
$$u-c'c^{-1}z'=ry^{m_0}+l''=ry'-rl_2+l''\in I,$$
and since  $l''-rl_2\in I\cap J$, it follows that  $u\in
Rz'\oplus (I\cap J)$, a contradiction.

Therefore  $I=Rz'\oplus (I\cap J)$. By Lemma 2.6, $I\cap J$ is a
direct sum of at most $|\Lambda|$ simple  $R$-modules. It follows
that $I$ is a direct sum of at most $|\Lambda|+1$ cyclic
$R$-modules by Lemma 2.7. We need only consider the structure of
each ideal of $R/{\rm Ann}(z')$. It is easy to check that $ {\rm
Ann}(z')= {\rm Ann}(x^{n_{0}-1})\cap  {\rm Ann}(y^{m_{0}-1})$.
Also, a trivial verification shows that
$${\rm Ann}(x^{n_{0}-1})=H_1\oplus Ry\oplus L~~~{\rm and}~~~
{\rm Ann}(y^{m_{0}-1})=Rx\oplus H_2\oplus L$$
 where\\
   $~~~~~~~~~~~~~~~~~H_1=\{rx\in
Rx~|~ rx+sy+l\in  {\rm Ann}(x^{n_{0}-1}), {\rm~for~some }~s\in
R,~l\in L \}$\\   $~~~~~~~~~~~~~~~~~H_2=\{sy\in Ry~|~ rx+sy+l\in
{\rm Ann}(y^{m_{0}-1}), {\rm~for~some }~r\in R,~l\in L \}.$

 Therefore  we conclude that $ {\rm Ann}(z')=H_1\oplus H_2\oplus
 L$. Put $\bar{R}=R/ {\rm Ann}(z')$. Then  $\bar{R}$ is a local ring
 with the maximal ideal   $\bar{\cal{M}}:={\cal{M}}/ {\rm Ann}(z')\cong Rx/H_1\oplus Ry/H_2$.
  It follows that  $\bar{\cal{M}}$  is a
   direct sum of at most two cyclic $\bar{R}$-modules. Thus, by
    Corollary 3.2, the ring $\bar{R}$ is Noetherian. Therefore, by
    \cite[Theorem 2.11]{Behboodi1}, every ideal of $\bar{R}$
     is a direct sum of at most two cyclic modules (see Result
     1.2)  which completes the proof.$~\square$

 Now we are in a
position to prove the main theorem of this paper. In fact we
  describe  the ideal structure of local rings $R$ for
which every ideal is a direct sum of cyclic modules.

\noindent{\bf Theorem  3.7.} {\it Let $(R, \cal{M})$  be a
 local  ring. The following
statements are equivalent:}\vspace{2mm}\\
(1) {\it Every ideal of $R$ is a direct sum of cyclic $R$-modules.}\vspace{1mm}\\
(2) {\it Every ideal of $R$ is  a direct sum of cyclic
$R$-modules, at most two of which
are not \indent simple.}\vspace{1mm}\\
(3) {\it There is an index set $\Lambda$ and a set of elements
$\{x, y\} \cup \{w_{\lambda}\}_{\lambda\in \Lambda}\subseteq R$
such that \indent ${\cal{M}}=  Rx  \oplus  Ry
\oplus(\bigoplus_{\lambda\in \Lambda} Rw_{\lambda})$ with: each
$Rw_{\lambda}$ a simple $R$-module and  $R/{\rm Ann}(x)$, \indent
$R/
{\rm Ann}(y)$  principal ideal rings.}\vspace{1mm}\\
(4) {\it There is an index set $\Lambda$ such that every  ideal
 $I$ of $R$ is one  of the following  forms:}\vspace{1mm}\\
\indent (i) {\it  $I=Rx'\oplus Ry'\oplus(\bigoplus_{\gamma\in
\Gamma} Rw'_{\gamma})$   where $\Gamma$ is an index set with
$|\Gamma|\leq|\Lambda|$ where \indent\indent  ${Rw'_{\gamma}}$'s
are
 simple $R$-modules and $ x'$, $y'\in R$ such  that
$R/{\rm Ann}(x')$,
$R/{\rm Ann}(y')$ \indent\indent are   principal ideal  rings.}\\
\indent (ii) {\it $I=Rz'\oplus(\bigoplus_{\gamma\in \Gamma}
Rw'_{\gamma})$ where $\Gamma$ is an index set with
$|\Gamma|\leq|\Lambda|$, ${Rw'_{\gamma}}$'s are simple
\indent\indent $R$-modules and $z'\in R$ is such that every ideal
of $R/{\rm Ann}(z')$ is a direct sum of at
 \indent\indent most  two principal ideals.}\vspace{1mm}\\
(5) {\it Every ideal of $R$ is a direct summand of a direct sum of
cyclic $R$-modules.}\vspace{1mm}\\ ({\it We note that the index
sets $\Lambda$'s   in the above statements {\rm (3)} and {\rm (4)}
are the same and  $|\Lambda|+2$ is a bound for the direct sum
decompositions of all ideals of $R$}).\vspace{1mm}\\

\noindent{\bf Proof.} $(1)\Rightarrow (3)$. Follows from Theorem 3.1 and Theorem 3.3.\\
$(3)\Rightarrow(4)$. Follows from Propositions 3.5 and 3.6.\\
 $(4)\Rightarrow (2)$, $(2)\Rightarrow (1)$ and $(1)\Rightarrow (5)$ are clear.\\
              $(5) \Rightarrow (1)$. Follows from   Lemma 2.5.$~\square$

By Corollary 3.2, Theorem 3.7 and \cite[Theorem 2.11]{Behboodi1},
we have the following ideal structure description for local rings
$(R, \cal{M})$ where $\cal{M}$ is finitely generated and every
ideal of $R$ is a direct sum of cyclic modules. Also, this result
is an analogue of Kaplansky's theorem \cite[Theorem
12.3]{Kaplansky1}.

\noindent{\bf Corollary  3.8.} {\it Let $(R, \cal{M})$  be a
 local  ring. Then the following
statements are equivalent:}\vspace{2mm}\\
(1) {\it $\cal{M}$ is finitely generated and every ideal of $R$ is a direct sum of cyclic $R$-modules.}\vspace{1mm}\\
(2) {\it $R$ is   Noetherian and every ideal
of $R$ is a direct sum of cyclic $R$-modules.}\vspace{1mm}\\
(3) {\it There exist  an positive integer  $n$ and a set of
elements $\{v_1,\cdots ,v_n\}\subseteq R$ such that
 \indent ${\cal{M}}=Rv_1\oplus\cdots\oplus Rv_n$  with: each $R/{\rm
Ann}(v_i)$ a principal ideal ring and at most
two  of \indent  $Rv_i$'s   not   simple.}\vspace{1mm}\\
(4) {\it  There exists  an positive integer  $n$  such that every
ideal of $R$ is of the form  $I=\indent Rv_1\oplus\cdots\oplus
Rv_m$
 where
 $m\leq n$, $\{v_1,\cdots,v_m\}\subseteq R$ and
 at most two of   $Rv_i$'s are \indent  not simple.}\vspace{1mm}\\
(5) {\it There exists  an positive integer  $n$  such that every
ideal of $R$ is a direct sum of at most \indent $n$
   cyclic   $R$-modules.}\vspace{1mm}\\
(6) {\it There exists  an positive integer  $n$ such that every
ideal of $R$ is a direct summand of a  \indent direct sum of  at
most $n$
cyclic $R$-modules.}\vspace{1mm}\\
 ({\it We note that the integers $n$'s in the above statements
are  the same and it is a bound for the direct sum  decompositions
of
all ideals of $R$}).\vspace{1mm}\\

 \noindent{\bf Remark  3.9.}  Let $R_1,\ldots,R_k$, where
 $k\in\Bbb{N}$,
  be  nonzero rings, and let
   $R$ denote the  direct product ring
   $\prod_{i=1}^kR_i$.  It is well-known that, if $I_i$ is an ideal of $R_i$ for
each $i=1,\ldots,k$, then $I=\prod_{i=1}^kI_i$ is an ideal of
 $R$. Furthermore  each ideal of $R$ has this form. It is straightforward
   to check that the ideal $I=\prod_{i=1}^kI_i$ of $R$  is a direct sum of cyclic ideals of $R$
if and only if the ideal  $I_i$ is a direct sum of cyclic ideals
of $R_i$ for each $i=1,\ldots,k$.

 We are thus led to the following
strengthening of Theorem 3.7 and Remark 3.9.

\noindent{\bf Corollary  3.10.} {\it Let $R=\prod_{i=1}^kR_i$,
where $k\in\Bbb{N}$ and where each $R_i$ is a  local ring with
maximal ideal ${\cal{M}}_i$ $(1\leq i\leq k)$. The following
statements are equivalent:}\vspace{2mm}\\
(1) {\it Every ideal of $R$ is a direct sum of cyclic $R$-modules.}\vspace{1mm}\\
(2) {\it For each $i$, every ideal of $R_i$ is  a direct sum of
cyclic $R_i$-modules, at most two of which
\indent are not simple.}\vspace{1mm}\\
(3) {\it For each $i$, there is an index set $\Lambda_{i}$ and a
set of elements $\{x, y\} \cup \{w_{\lambda}\}_{\lambda\in
\Lambda_{i}}\subseteq R_i$ \indent such  that
${\cal{M}}_i=R_ix\oplus R_iy\oplus(\bigoplus_{\lambda\in
\Lambda_i} R_iw_{\lambda})$ with: each $R_{i}w_{\lambda}$ a simple
$R$-module and \indent $R_{i}/{\rm Ann}(x)$,
$R_{i}/ {\rm Ann}(y)$  principal ideal rings.}\vspace{1mm}\\
(4) {\it There exit index sets $\Lambda_1,\ldots, \Lambda_k$ such
that for each $1\leq i\leq k$, every ideal of $R_i$ is  of \indent
the forms  $I=R_ix'\oplus R_iy'\oplus(\bigoplus_{\gamma\in
\Gamma_i} R_iw'_{\gamma})$ or
 $I=R_iz'\oplus(\bigoplus_{\gamma\in \Gamma_i}
R_iw'_{\gamma})$ where $\Gamma_i$ is an  \indent index   set with
$|\Gamma_i|\leq|\Lambda_i|$, $R_iw'_{\gamma}$'s (${\gamma\in
\Gamma_i}$)  are  simple $R_i$-modules,  $x'$, $y'$, $z'\in R_i$
 such \indent that  $R_i/{\rm Ann}(x')$,   $R_i/{\rm Ann}(y')$ are principal
ideal rings and  every ideal of $R_i/{\rm Ann}(z')$ \indent is a
direct sum of at
  most  two principal ideals.}\vspace{1mm}\\
(5) {\it Every ideal of $R$ is a direct summand of a direct sum of
cyclic $R$-modules.}\\

\noindent{\bf 4. Examples}

\noindent In this section we provide several examples illustrating
the results of Section 3 as well as the necessity of certain
hypotheses in these results. We begin  with the following
interesting example. In fact, the following example shows that the
corresponding of the above Corollary 3.10, is not true in general
for  the case $R=\prod_{\lambda\in\Lambda}R_\lambda$ where
$\Lambda$ is an infinite index set and each $R_\lambda$ is a local
 ring (even if for each $\lambda\in\Lambda$,
 $R_\lambda\cong\Bbb{Z}_2$).

  \noindent{\bf Example 4.1.}  Let $\Lambda$ be an infinite index set
  and  $\{F_{\lambda}~|~{\lambda\in\Lambda}\}$ be a set of fields.  We put
    $R=\prod_{\lambda\in\Lambda}F_\lambda$.  Thus for each
   ${\lambda\in\Lambda}$, ${F_\lambda}$ is a local ring and
   every ideal of ${F_\lambda}$ is   cyclic (that is $(0)$ or
   ${F_\lambda}$). Clearly, $I=\bigoplus_{\lambda\in\Lambda}F_\lambda$  is a non-maximal  ideal of $R$ and hence
    there exists a
   maximal ideal $P$ of $R$ such that $I\subsetneq P$.  It was shown by Cohen and Kaplansky
    \cite[Lemma 1]{Cohen} that $P$ is not a direct sum of
   principal ideals. Thus  the
corresponding of the above Corollary 3.10, is not true in general
for  the case $R=\prod_{\lambda\in\Lambda}R_\lambda$ where
$\Lambda$ is an infinite index set.

Let $(R, \cal{M})$ be a local ring.  By Theorem 3.7,  every ideal
of $R$ is a direct sum of cyclic $R$-modules if and only if there
is an index set $\Lambda$ and a set of elements $\{x, y\} \cup
\{w_{\lambda}\}_{\lambda\in \Lambda}\subseteq R$ such that
${\cal{M}}=Rx  \oplus  Ry \oplus(\bigoplus_{\lambda\in \Lambda}
Rw_{\lambda})$ with: (i) each $Rw_{\lambda}$ a simple $R$-module,
 (ii) $R/  {\rm Ann}(x)$, $R/ {\rm Ann}(y)$  principal ideal rings.
Moreover, in this case,  every ideal $I$ of $R$ has the  form
$I=Rx'\oplus Ry'\oplus(\bigoplus_{\gamma\in \Gamma} Rw'_{\gamma})$
where $\Gamma$ is an index set with $|\Gamma|\leq|\Lambda|$,
${Rw'_{\gamma}}$s ($\gamma\in \Gamma$) are  simple $R$-modules and
$x'$, $y'\in R$. The following example shows that  property (ii)
does not hold for all ideals of $R$, even if $R$ is Artinian and
${\cal{M}}$ is two generated.

\noindent{\bf Example 4.2.} Let $F$ be a field, let $n\geq 3$  and
let  $R$ be the $F$-algebra with generators $x$, $y$ subject to
the relations $x^n=y^n=xy=0$, $R\cong F[X,Y]/<X^n,Y^n,XY>$. The
ring $R$ is a Noetherian local ring with maximal ideal
${\cal{M}}=Rx\oplus Ry$.  Since ${\cal{M}}^n=(0)$,  dim$(R)=0$ and
so $R$ is an Artinian local ring.  By Theorem 3.7, every ideal of
$R$ is a direct sum of at most two cyclic $R$-modules and  $R/
{\rm Ann}(x)$, and $R/ {\rm Ann}(y)$ are principal ideal rings.
Let $I=Rz$ where $z=x+y$. We note that  if
$I=\bigoplus_{i=1}^nRz_i$ where $n\in \Bbb{N}$ and
 $Rz_i$ are nonzero cyclic  $R$-modules,   then   by Lemma 2.7,
   $n=1$. Set $\bar{R}=R/{\rm Ann}(z)$. Clearly
$ {\rm Ann}(z)=Rx^{n-1}\oplus Ry^{n-1}$. Since
$\bar{\cal{M}}:={\cal{M}}/{\rm Ann}(z)\cong (Rx/Rx^{n-1})\oplus
(Ry/Ry^{n-1})$, it follows that the maximal ideal  $\bar{\cal{M}}$
of $\bar{R}$ is a    direct sum of  two nonzero cyclic
$\bar{R}$-modules and hence  by Lemma 2.7, ${\cal{M}}$ is not
principal, i.e., $\bar{R}$ is not a principal ideal ring.

 By \cite[Example 3.1]{Behboodi1}, for each integer $n\geq 3$, there exists  an
Artinian (Noetherian) local ring $(R, {\cal{M}})$ such that
${\cal{M}}$ is a direct sum of $n$ cyclic $R$-modules, but there
exists a two generated ideal of $R$ such that it  is not a direct
sum of cyclic $R$-modules. Next, the following example shows that
there exists also a non-Noetherian local ring $R$ such that every
prime ideal of $R$ is a direct sum of cyclic $R$-modules, but some
of the ideals of $R$ are  not direct sums of cyclic $R$-modules.

 \noindent{\bf Example 4.3.} Let $F$ be a field and let $R$ be
 the $F$-algebra with generators $\{x_i~|~i\in\Bbb{N}\}$ subject to the following relations
$x_1^3=x_2^3=x_3^3=x_k^2=0,~~k\geq 4~ {\rm and}~  x_ix_j=0~{\rm
for~all}~ i\neq j,$  $R\cong F[\{X_i~|~i\in\Bbb{N}\}]/<X_1^3,
X_2^3, X_3^3, X_k^2, X_iX_j~|~4\leq k\in\Bbb{N}, ~ i\neq j\geq
1>$. Then $R$ is a non-Noetherian local ring with the maximal
 ideal ${\cal{M}}=\bigoplus_{i\in \Bbb{N}} Rx_i$. Clearly Spec$(R)=\{{\cal{M}}\}$ since ${\cal{M}}^3=(0)$.
   Thus  every   prime ideal of $R$
   is a  direct sum of  cyclic   $R$-modules, but by Lemma 2.3, the ideal $J=R(x_1+x_2)+R(x_1+x_3)$
  is not  a direct sum of cyclic $R$-modules.

By Theorem 3.1, if a local ring $(R, \cal{M})$ has the property
that every   ideal $I$ of $R$  is a direct sum of cyclic
$R$-modules, then dim$(R)\leq 1$ and  $|{\rm{Spec}}(R)|\leq 3$.
Clearly, $|{\rm{Spec}}(R)|=1$ whenever  dim$(R)=0$ and
$|{\rm{Spec}}(R)|=2$ or $3$  whenever dim$(R)=1$.  The following
examples cover all the different cases mentioned above for
dim$(R)$ and $|{\rm{Spec}}(R)|$.

 \noindent{\bf Example 4.4.}  Let $F$ be a field and $n\in\Bbb{N}$. Consider the following
 rings:\vspace{2mm}\\
 (1) $R_1=F[[X]]$ (formal power series ring).\\
  (2) $R_2=F[\{X_i~|~1\leq i\leq n\}]/<\{X_iX_j~|~1\leq i,j\leq
 n\}>$.\\
 (3) $R_3=F[\{X_i~|~i\in\Bbb{N}\}]/<\{X_iX_j~|~i,j\in
 \Bbb{N}\}>$.\\
 (4) $R_4=F[[X,Y]]/<XY, Y^2>$.\\
 (5) $R_5=F[[\{X_i~|~i\in\Bbb{N}\}]]/<\{X_iX_j~|~i\neq ~j\}\cup \{X_i^2~|~i\geq 2\}>$.\\
(6) $R_6=F[[X,Y]]/<XY>$.\\
(7) $R_7=F[[\{X_i~|~i\in\Bbb{N}\}]]/<\{X_iX_j~|~i\neq ~j\}\cup
\{X_i^2~|~i\geq 3\}>$.

   It is  easy to check that all above
rings are local and  by Theorem 3.7, one can easily see that for
each $i\in \{1,2,\cdots,7\}$ the  ring $R_i$ has the property that
every ideal is a direct sum of cyclic ideals. Let ${\cal{M}}_i$
denote the maximal ideal of $R_i$ for $i\in \{1,2,\ldots,7\}$.
 Then  we easily obtain the following:

  \noindent (1) $R_1$ is a  domain (in fact $R_1$ is a local  PID) with  dim$(R_1)=1$,
  ${\cal{M}}_1=<X>$  and
   \indent ${\rm{Spec}}(R_1)=\{(0), {\cal{M}}_1\}$.

  \noindent (2) $R_2$ is a non-domain   Artinian ring with  dim$(R_2)=0$,
  ${\cal{M}}_2=R_2x_1\oplus\cdots \oplus R_2x_n$
 (where \indent $x_i=X_i+<\{X_iX_j~|~1\leq i,j\leq
 n\}>$) and  ${\rm{Spec}}(R_2)=\{{\cal{M}}_2\}$.

  \noindent (3) $R_3$ is a non-domain,    non-Noetherian  ring with dim$(R_3)=0$,
   ${\cal{M}}_3=\bigoplus_{i=1}^\infty R_3x_i$ (where \indent
 $x_i=X_i+<\{X_iX_j~|~i,j\in  \Bbb{N}\}>$)   and
 ${\rm{Spec}}(R_3)=\{{\cal{M}}_3\}$.

\noindent (4) $R_4$ is a non-domain,   Noetherian ring with
dim$(R_4)=1$, ${\cal{M}}_4=R_4x\oplus R_4y$ (where \indent
$x=X+<XY, Y^2>$ and $y=X+<XY, Y^2>$) and
${\rm{Spec}}(R_4)=\{{\cal{M}}_4, R_4y\}$.

\noindent (5) $R_5$ is a non-domain,  non-Noetherian ring with
dim$(R_5)=1$,
 ${\cal{M}}_5=\bigoplus_{i=1}^\infty R_5x_i$ (where
 \indent $x_i=X_i+<\{X_iX_j~|~i\neq ~j\}\cup \{X_i^2~|~i\geq 2\}>$) and
 ${\rm{Spec}}(R_5)=\{{\cal{M}}_5, \bigoplus_{i=2}^\infty R_5x_i\}$.

  \noindent (6) $R_6$ is a non-domain,   Noetherian ring with  dim$(R_6)=1$,  ${\cal{M}}_6=R_6x\oplus R_6y$
  (where \indent $x=X+<XY>$ and $y=Y+<XY>$)
  and ${\rm{Spec}}(R_6)=\{{\cal{M}}_6,~R_6x,~  R_6y\}$.

 \noindent (7) $R_7$ is a non-domain,  non-Noetherian ring with dim$(R_7)=1$,
${\cal{M}}_7=\bigoplus_{i=1}^\infty R_7x_i$ (where \indent
$x_i=X_i+<\{X_iX_j~|~i\neq ~j\}\cup \{X_i^2~|~i\geq 3\}>$) and
   \begin{center}${\rm{Spec}}(R_7)=\{{\cal{M}}_7,~ \bigoplus_{i=2}^\infty R_7x_i,~R_7x_1\oplus (\bigoplus_{i=3}^\infty
   R_7x_i)\}.$\end{center}

\end{document}